\documentclass[a4paper,10pt]{article}
\def \Z {{\mathbf {Z}}}
\def \R {{\mathbf {R}}}
\def \N {{\mathbf {N}}}

\title{\bf  Абсолютная непрерывность  и сингулярность  \\ спектра
 потоков $\bf T_t\otimes T_{at}$}
\author{В.В.Рыжиков}
\date{}
\usepackage[T1]{fontenc} 
\usepackage[cp1251]{inputenc}  
\usepackage[russian]{babel}
\usepackage{graphicx}
\graphicspath{{pictures/}}
\DeclareGraphicsExtensions{.pdf,.png,.jpg}
\begin{document}
\large
\maketitle
\begin{abstract}  

\vspace{2mm}
Отвечая на вопрос В.И. Оселедца, мы предъявляем случайную величину $\xi$ такую,
что суммы вида  $\xi(x)+a\xi(y)$ имеют  сингулярное распределение
для всюду плотного в $(1, +\infty)$ множества  значений параметра $a$, причем для другого
всюду плотного  множества соответствующие распределения абсолютно непрерывны. 
Предлагаемые  распределения имеют  динамическое происхождение. 
Они реализуются   спектральными мерами   подходящих эргодических потоков на пространстве с сигма-конечной мерой.

Библиография: 5 названий, УДК: 517.987, \ MSC: Primary 28Y05; Secondary 58F11

Ключевые слова и фразы: \it распределение случайной величины, спектр тензорных произведений потоков,
диссипативость, слабые пределы операторов.\rm

\end{abstract} 

\section{Введение}

 В связи с изучением мер Эрдёша (см. \cite{O}, \cite{B}) В.И. Оселедец предложил  задачу, которую 
 мы формулируем следующим образом.
 
\it Пусть заданы  $C$ и $D$ -- непересекающиеся   счетные плотные подмножества луча $(1, +\infty)$. 
Требуется  найти меру $\sigma$   на прямой $\R$ такую, что 
 произведение  $\sigma\times \sigma$ проектируется в сингулярные меры   на горизонтальной прямой 
вдоль  направлений $(c,1)$, $c\in C$, и проектируется в абсолютно непрерывные   меры  вдоль 
 направлений $(d,1)$,  $d\in D$.\rm

В заметке  показано, что нужными свойствами обладают   спектральные меры некоторых потоков  
на пространстве с сигма-конечной мерой.

\vspace{3mm}
\bf Теорема.  \it Пусть $C,D$ -- счетные непересекающиеся плотные подмножества луча $(1,+\infty)$.
Найдется сохраняющий сигма-конечную меру поток $T_t$ такой, что  
автоморфизмы $T_1\otimes T_{c}$ имеют простой сингулярный спектр для каждого $c\in C$,
а автоморфизмы $T_1\otimes T_{d}$ имеют счетнократный лебеговский  спектр для всех $d\in D$.
\rm

\vspace{3mm}
Напомним, что унитарный поток $T_t$, обладающий циклическим вектором, 
для некоторой конечной  борелевской меры $\sigma$ на $\R$  изоморфен потоку  $U_t$, 
определенному формулой
 $$ U_t f(s) = e^{its}f(s), \ \ f\in L_2(\R, \sigma).$$ 
Мера $\sigma$ (и всякая  эквивалентная ей мера)  называется спектральной мерой потока $T_t$.

  Обеспечить    сингулярность спектра $T_t\otimes T_{ct}$ при $c\in C$ удобно, используя 
слабые пределы. Если  найдется последовательность ${t_j}$, для которой
$$T_{t_j}\otimes T_{ct_j}\to_w \frac  {I\otimes I} 4,  \eqno (1)$$
 то спектр  потока $T_t\otimes T_{ct}$ будет сингулярным.  Иначе для некоторого вектора $F\neq 0$ мы получили бы  
$(T_{t_j}\otimes T_{ct_j})F\to_w 0, $
что невозможно. 

При $d\in D$  для некоторого множества $Y$ положительной меры
и числа $N(d)$  для эргодического потока $T_{t}$  будет обеспечено условие
$$\forall  t>N(d)\ \ \ (Y\times Y)\cap  (T_{t}Y\times T_{d t}Y)\ =\ \phi. \eqno (2)$$
Это  влечет за собой диссипативность потока $T_t\otimes T_{at}$ и, как следствие,
 абсолютную непрерывность его спектра. Напомним, что диссипативность потока $S_t$ означает
существование такого измеримого множества $W$, что множества $S_nW$, $n\in\Z$, не пересекаются,
а их объединение есть все пространство, на котором действует поток.

Спектральная мера  потока с отмеченными выше свойствами является решением задачи Оселедца.
В  \cite{R} приведены примеры автоморфизмов $T$, для которых 
сингулярность и лебеговость  спектра произведения $T\otimes T^n$ зависела от $n\in \N$.
Ниже предлагается модификация этих примеров.

\section{Конструкция потока}
Поиск подходящих примеров происходит в   классе  потоков ранга один,   заданных   параметрами
$h_1, w_1>0$  и последовательностями  неотрицательных чисел  $s_j(1)$, $s_j(2)$, $s_j(3)$, $s_j(4)$,
 $j\in\N$.  Определим   индуктивно фазовое пространство и поток, отвечающие указанным параметрам. 
 На шаге  1 дан прямоугольник $X_1$ ширины  $w_1$ и высоты  $h_1$.
Считаем, что  $X_1$ (и все так называемые прямоугольники  $X_j$, определенные  ниже) 
 не содержат правой границы. 
На $X_1$  поток задан как движение   вверх по вертикали 
с постоянной единичной скоростью. Движение точки после достижения 
верхней границы определяется далее.
Пусть на шаге $j$ построена  башня $X_j$, которая отождествляется с прямоугольником
ширины $w_j={w_1}{4^{1-j}}$  и высоты
$
h_j = 4 h_{j-1} + \sum_{i=1}^{4} s_j (i).
$ 
На шаге $j+1$ прямоугольник $X_j$ разрезается по вертикали на 
$4$ одинаковых прямоугольника $X_j^i$, которые называются колоннами.  Над колоннами $X_j^i$, $i=1,2,3,4$,  
надстраиваем прямоугольники  высоты  $s_j (i)$, соответственно.
Получаем $4$ новых прямоугольника $\tilde X_j^i$ (надстроенные колонны) ширины ${w_1}{4^{-j}}$.
При $i<4$ склеиваем  верхнюю границу колонны $\tilde X_j^i$  с нижней границей 
колонны  $\tilde X_j^{i+1}$, $i=1,2,3$. 
Объединение  этих колонн по определению есть   башня $X_{j+1}$,
 которую мы ассоциируем  с прямоугольником высоты $h_{j+1}$ и ширины $w_{j+1}={w_1}{4^{-j}}$. 
Поток в этом прямоугольнике  определен как движение  вверх по вертикали с постоянной единичной скоростью.

На самом деле, $X_j$ -- это часть плоскости $\bf{R}^2$, которую (для удобства определения 
и  изучения свойств полученного специального потока) отождествляют с прямоугольником. 
Фазовое пространство $X$ потока $T_t$  -- это объединение 
$\cup_{j=1}^{\infty} X_j$. Движение точки под действием потока определено 
в башне $X_j$ до тех пор, пока точка не покидает башню  $X_j$. Потом рассматриваем ее движение
 в  башнях $X_{j+1}$, $X_{j+2}$ и т. д.   По построению для каждого  $t$ преобразование 
$T_t$ является биекцией на $X$, сохраняющей плоскую меру Лебега $\mu$. Про потоки ранга один известно,
что они эргодичны (нет нетривиальных инвариантных измеримых подмножеств) и обладают  простым спектром 
(есть циклический ветор).  Этими фактами мы  воспользуемся позже.

Пусть параметры конструкции удовлетворяют следующим условиям: 
$$ s_j(1)=0, \ \ s_j(3)=(c_j-1)h_j,  \ \ s_j(4)\gg s_j(2)\gg h_j. $$
Здесь выражение вида $s_j(2)\gg h_j$ означает, что отношение $s_j(2)/ h_j\to\infty$ 
при $j\to\infty$. Мера фазового пространства такого потока бесконечна.
Последовательность $c_j$ принимает все значения из множества $C$, причем каждое значение 
она принимает бесконечное  число раз. Условия $ s_j(1)=0, s_j(3)=(c_j-1)h_j$    обеспечат свойство (1).    Выбор быстро растущих последовательностей $s_j(2)$, $s_j(2)$  при
условии быстрого роста отношения   $s_j(4)/ s_j(2)$ позволит реализовать свойство (2) для всех $d\in D$.

\section{Доказательство теоремы}
\bf Сингулярность спектра.   \rm Начнем с проверки того, что при
 $c\in C$ автоморфизм $T_1\otimes T_{c}$ обладает свойством: 
 для измеримых множеств $A,B$, которые являются объединениями 
прямоугольников ширины $w_k$ в башне $X_k$ ($k$ фиксировано) и всех достаточно больших  $j$ верно 
$$\mu (T_{h_j}A\cap B)= \mu (A\cap B)/4, \ \ \mu (T_{c_jh_j}A\cap B)= \mu (A\cap B)/4. \eqno (3)$$
Множество  $A$ является объединением множеств $A^i=X_j^i\cap A$, $i=1,2,3,4$.
Из $s_j(2) / h_j\to\infty$, $c>1$ и $s_j(4) / h_j\to\infty$ получаем, что образы $T_{h_j}A^i$ при $i=2,3,4$ попадают в надстройки над $X_j$, поэтому они не пересекаются с $B\subset X_j$. 
Таким  образом,   имея
$$T_{h_j}A \cap B=T_{h_j}A^1 \cap B= A\cap B\cap X_j^2,$$
получаем первое равенство из (3).
Аналогично для  больших $j$ устанавливается второе равенство из (3).
Так как линейные комбинации индикаторов рассмотренных множеств $A,B$ плотны в $L_2(\mu)$, 
а последовательность $c_j$ такова, что для каждого $c\in C$ найдется последовательность $j_k$,
для которой $c_{j_k}=c$, получаем 
$$(T_1\otimes T_{c})^{h_{j_k}}\to_w \ \frac  {I\otimes I}{16}.$$ 
Спектр автоморфизма  $T_1\otimes T_{c}$ сингулярен, что 
влечет за собой сингулярность спектра потока  $T_t\otimes T_{ct}$.

\vspace{2mm}
\bf Абсолютная непрерывность  спектра.   \rm Теперь  надо показать, что при подходящем выборе параметров $s_j(2), s_j(4)$
для всех $d\in D$ поток $T_t\otimes T_{dt}$ диссипативен.
На шаге $j$ мы выбираем $s_j(2), s_j(5)$  так, чтобы  для  $Y=X_1$ выполнялось
$$\forall \ k\leq j \ \forall \ t\in [h_j, h_{j+1}] \ \ (T_{t}Y\times T_{d_k t}Y)\cap (Y\times Y)=\ \phi. $$

Так как $h_j\gg h_{j-1}$,
 множество $Y\subset X_{j-1}$ расположено для  больших $j$ вблизи  основания башни  $X_j$.
  Из этого вытекает, что условие 
$\mu (T_{t}Y\cap Y)> 0$, $h_j\leq t \leq h_{j+1}$ выполнено только в   случаях, когда $t\in [h_j,h_{j+1}]$
 близко к одному из чисел: $h_j$, $c_j h_j$,$s_j(2)$, $s_j(4)$  (мы называем два числа близкими,
если их отношение близко к 1). 

Условие 
$\mu (T_{d_kt}Y\cap Y)> 0$, $h_j\leq t \leq h_{j+1}$ по аналогичной причине 
выполнено только в   случаях, когда $t$
относительно близко к одному из чисел: $d_kh_j$, $d_kc_j h_j$,  $d_ks_j(2)$, $d_ks_j(4)$.
При выборе достаточно  больших   $s_j(2), s_j(4)$  и  $\frac {s_j(4)}{s_j(2}$, 
с учетом того, что набор  чисел $d_1, d_2,\dots, d_j$ не содержит 1 и $c_j$,
получаем желаемое:  числа из набора  $h_j$, $c_j h_j$, $s_j(2)$, $s_j(4)$ настолько отличаются от чисел
 набора $d_kh_j$, $d_kc_j h_j$, $d_ks_j(2)$, $d_ks_j(4)$,
что условие $\mu (T_{t}Y\cap Y)> 0$ несовместимо с условием $\mu (T_{d_kt}Y\cap Y)> 0$
при $k=1,2,\dots,j$,  $t\in [h_j,h_{j+1}]$. 

Таким образом, мы построили  поток $T_t$ такой, что для $d=d_k\in D$  для всех $n> h_k$ выполнено 
$$(T_{n}Y\times T_{d n}Y)\cap (Y\times Y)=\ \phi.$$
В силу этого свойства  множества $Y\times Y$, эргодичности потока $T_{t}$ и того, что автоморфизм 
$T_{1}\times T_{d}$ коммутирует с $Id\times T_t$ при $t\in\R$, вытекает   диссипативность автоморфизма $T_{1}\times T_{d}$.  Его спектр счетнократный лебеговский,
а это гарантирует  абсолютную непрерывность спектральной меры потока $T_{t}\otimes T_{d t}$
(на самом деле спектр такого  потока счетнократный лебеговский).

\vspace{2mm}
\bf Простота сингулярных спектров. \rm 
Чтобы добиться свойства простоты спектра, нам понадобится небольшая модификация описанной  ранее конструкции потока.
Рассмотрим следующие  возмущения параметров $s_j(1),s_j(3)$, положив
$$s_j(1)=\Delta_j(1), \ \   s_j(3)=(c_j-1)h_j + \Delta_j(3).$$ 
При этом неотрицательные последовательности  $\Delta_j(1), \Delta_j(3)\leq 1$ выбираются так, чтобы для любого $c\in C$ замыкание множества 
$\{(\Delta_{j_k}(1), \Delta_{j_k}(3)): b_{j_k}=c\}$ было квадратом   $[0,1]^2$. 
Тогда  для  $a,b\in [0,1]$ найдется последовательность $j'$ такая, что  
$\Delta_{j'}(1)\to a$, $\Delta_{j'}(3)\to b$, $c_{j'}=c$, а это приводит к 
$$(T_1\otimes T_{c})^{h_{j'}}\to_w \frac 1 {16}T_{-a}\otimes T_{-b}.$$
При условии, что спектр потока $T_t$ простой, такие слабые   пределы  обеспечивают  простоту 
спектра оператора  $T_1\otimes T_{c}$. 
Покажем это.  Алгебра фон Неймана, порожденная оператором   $T_1\otimes T_{c}$, в силу сказанного будет содержать
все операторы вида $T_a\otimes T_b$ при $a,b\in [-1,0]$. Она замкнута   относительно операций сопряжения и умножения,  следовательно, будет содержать все операторы вида $T_s\otimes T_t$, $s,t\in \R$. 
Пусть $f$ -- циклический вектор для потока $T_t$.  Тогда $f\otimes f$  -- 
циклический вектор для оператора $T_1\otimes T_{c}$.  Действительно,  линейное и топологическое замыкание множества 
$\{T_n f\otimes T_{cn}f : \ n\in\N\}$  содержит все векторы вида 
$T_sf\otimes T_tf$, $s,t\in \R$, поэтому оно является пространством  $L_2(\mu)\otimes L_2(\mu)$.

При  ограниченном возмущении параметров  $s_j(1),s_j(3)$ свойство  диссипативности автоморфизмов  
$T_1\otimes T_{d}$, $d\in D$, очевидно,  сохраняется.    Теорема доказана.

\newpage
\bf Замечания.  \rm Если множество $C$ конечно, то в рамках нашей конструкции автоморфизмы 
 $T_1\otimes T_{d}$ диссипативны для всех 
$d\in (1, +\infty)\setminus C$. Если же  $C$ плотно в $(1, +\infty)$, то  
   возникает дополнительное плотное $G_\delta$-множество параметров, отвечающих простому сингулярному спектру. 
При  этом множество параметров, отвечающих абсолютно непрерывному спектру, может иметь полную меру.
Тем самым возникает ситуация, когда для параметрического семейства преобразований  $T_1\otimes T_{a}$
сингулярность спектра  типична по Бэру, а лебеговость спектра   типична  по Лебегу.

  Отметим, что в работе \cite{LR} реализация  простоты сингулярного спектра 
$T_1\otimes T_{a}$  для всех $a\in (1,+\infty)$  представляла собой нетривиальную задачу, а
ее решение   также  использовало   технику слабых пределов.

\vspace{3mm}
 Автор благодарит Бенжамина  Вейса,  сообщившего в ходе обсуждения задачи Оселедца идею другого решения (совместно
с М. Хохманом), использующего  результаты работы \cite{N}.
 



\end{document}